\documentclass[journal]{IEEEtran}
\usepackage{cite}
\usepackage{url}
\usepackage{empheq}
\usepackage{float}
\usepackage{dsfont}
\usepackage{amsthm}
\usepackage{algorithm}
\usepackage{algorithmic}

\usepackage{amsmath,amssymb,amsfonts}
\usepackage{graphicx}
\usepackage{textcomp}
\usepackage{amsmath}
\usepackage{enumerate}
\usepackage{epstopdf}
\usepackage{array}
\usepackage{booktabs}
\usepackage{subfigure}
\usepackage{multirow}
\usepackage{soul, color}
\usepackage[usenames,dvipsnames]{xcolor}
\usepackage[latin1]{inputenc}

\newtheorem{proposition}{Proposition}

\soulregister\cite7 
\soulregister\citep7 
\soulregister\citet7 
\soulregister\ref7 
\soulregister\pageref7 

\begin{document}

\title{A Distributed Online Algorithm for Promoting \\Energy Sharing Between EV Charging Stations}

\author{
Dongxiang~Yan and
Yue~Chen

\thanks{D. Yan and Y. Chen are with the Department of Mechanical and Automation Engineering, the Chinese University of Hong Kong, Hong Kong SAR, China (e-mail: yuechen@mae.cuhk.edu.hk).
}}
\markboth{Journal of \LaTeX\ Class Files,~Vol.~XX, No.~X, Feb.~2019}%
{Shell \MakeLowercase{\textit{et al.}}: Bare Demo of IEEEtran.cls for IEEE Journals}
\maketitle

\begin{abstract}
In recent years, electric vehicle (EV) charging station has experienced an increasing supply-demand mismatch due to its fluctuating renewables and unpredictable charging demand. To reduce its operating cost, this paper proposes a distributed online algorithm to promote the energy sharing between charging stations. We begin with the offline and centralized version of the EV charging stations operation problem, whose objective is to minimize the long-term time-average total cost. Then, we develop an online implementation approach based on the Lyapunov optimization framework. Although the proposed online algorithm runs in a prediction-free manner, we prove that by properly choosing the parameters, the time-coupling constraints remain to be satisfied. We also provide a theoretical bound for the optimality gap between the offline and online optimums. Furthermore, an improved alternating direction method of multipliers (ADMM) algorithm with iteration truncation is proposed to enable distributed computation. The proposed algorithm can protect privacy while being suitable for online implementation. Case studies validate the effectiveness of the theoretical results. Comprehensive performance comparisons are carried out to demonstrate the advantages of the proposed method.

  \end{abstract}


\IEEEpeerreviewmaketitle

\vspace{-0.cm}
\section{Introduction}
\IEEEPARstart {D}{riven} by the increasing number of electric vehicles (EVs), massive EV charging stations need to be built to meet the growing charging demand~\cite{yang2019optimal}.
The charging stations can equip with renewable generation (e.g., photovoltaic (PV)) and battery energy storage to reduce the operating cost~\cite{yan2018optimized}.
The flexible EV charging demands and renewable generations allow charging stations to play an important role in the future energy systems.
Meanwhile, charging station changes to be an energy prosumer from a pure consumer.
The power supply-demand mismatch from intermittent renewable generation and stochastic EV charging load poses new challenges to its operation~\cite{wang2017stochastic}.
It is crucial to develop an effective method to mitigate the mismatch of supply and demand and improve its self-consumption.
Considering that different EV charging stations have distinct supply and demand patterns, energy sharing between them to make use of their complementary features has been considered as a promising solution~\cite{tushar2020peer}.

Recently, there have been extensive studies on energy sharing in the literature.
An energy trading strategy was developed for interconnected microgrids in~\cite{wang,2020fixLoad}.
The microgrid with surplus renewable energy can trade with others in need to gain mutual benefits.
An energy sharing framework designed for an energy building cluster was studied in~\cite{2018cuiBuilding} considering the building thermal dynamics.
A leader-follower game was used to model the interaction of an energy sharing provider and prosumers in a peer-to-peer (P2P) market~\cite{markov}.
The energy sharing model was extended to multi-energy systems in~\cite{jing2020fair}.
The generalized demand function based energy sharing scheme developed in~\cite{chen2020approaching} was proven to have a near socially optimal economic efficiency.
To encourage energy sharing between various prosumers, Nash bargaining based model~\cite{fan2018bargaining} and mid-market ratio mechanism~\cite{tushar2018peer} have been explored to allocate the payment among participants.
The above studies adopt deterministic approaches without uncertainties.
Some works considering uncertainties have been carried out.
For example, stochastic programming was used to treat the uncertain PV generation, load, and electricity price in the day-ahead energy sharing decision-making stage~\cite{stochastic}.
A two-stage robust model for energy sharing was developed~\cite{2019bi}.
However, the above works usually operate in an offline manner, adopt a predetermined time-of-use price, and assume complete information of the uncertainties. In practice, those data may be unavailable or inaccurate, and thus, the obtained solution may fail to adapt to the changes in the real-time.

Various online optimization based methods have been used for real-time implementation.
A common practice is the greedy algorithm which is shortsighted because it directly decomposes the time accumulated problem into each single time period ignoring the time-coupling constraints~\cite{guo2021asynchronous}.
Another method is the model predictive control (MPC).
An online MPC based optimal charging strategy was proposed for multiple charging stations to minimize the utility cost~\cite{zheng2018online}.
However, MPC still relies on a short-term forecast and the result is affected by the forecast accuracy~\cite{stai2020online}. Moreover, the rolling optimization can be computational expensive. 
An alternative real-time scheme is based on the Lyapunov optimization that takes the long-term benefit into account but requires no prior knowledge of uncertainty~\cite{fan2020online}.
Its decision-making just depends on the current state which is more flexible and practical.
The Lyapunov optimization based real-time scheme has been widely applied in many fields such as online network resource allocation~\cite{chen2017stochastic}, energy management in data center~\cite{yu2018distributed} and microgrid~\cite{2017RTMG}.
There are several works related to energy sharing.
For example, Lyapunov optimization was used to improve PV consumption of a cluster of nanogrids~\cite{2018Nano}. The energy trading between an end-user and grid was optimally scheduled to maximize the profit via Lyapunov optimization~\cite{Liu2020dynamic}.
However, the above works were based on the centralized scheme with concerns about privacy and communication burden. 
Distributed optimization methods, such as alternating direction method of multipliers (ADMM) \cite{2011boyd}, can help tackle these concerns. A real-time ADMM based algorithm was proposed to address the power balancing problem in a renewable-integrated power grid~\cite{sun2016distributed}.
Ref.~\cite{zhong2019online} applied ADMM to operate a shared energy storage system for multiple users in real-time.
An ADMM based online algorithm was developed to perform distributed energy management for multiple data centers~\cite{zhang2020distributed}.
However, the ADMM method may require a large number of iterations to converge. 
It would not be a big issue for the operation with a relatively long time scale. However, with the increasing penetration of renewable energy, the real-time operation runs in a smaller time resolution to cope with the volatility. Therefore, less and deterministic iterations are preferred, otherwise, the ADMM may not have enough time to converge.
Motivated by the above discussions, a research gap can be found: an EV charging stations energy sharing framework with online optimization to adapt to uncertainties, distributed implementation to protect privacy, and reduced iterations for fast execution has not been well explored yet.

In this paper, we propose a distributed online optimization for energy sharing between EV charging stations. Our main contributions are three-fold:
\begin{enumerate}
  \item We propose an online implementation approach to manage the energy sharing between EV charging stations. First, an offline optimization that minimizes the long-term time-average total cost is formulated. To fit into the Lyapunov optimization framework, two virtual queues are introduced to turn the inter-temporal battery energy and charging demand dynamics into mean rate stable conditions. Then, an augmented objective function with drift-plus-penalty term is built, based on which the online algorithm is derived. We prove that if the parameters are chosen within a particular range, the battery energy constraint is satisfied even if it is not explicitly considered. A theoretical bound for the optimality gap between the offline and online optimums is given. The proposed online algorithm requires no prior knowledge of the future, which is more practical.
  \item An improved ADMM based algorithm with iteration truncation is proposed. Distinct from existing work, our algorithm's iteration ends when it reaches a threshold, followed by a shared energy balance update mechanism. This can greatly cut down computational time to fit into online applications. In simulation, we find that compared with the conventional ADMM, the proposed algorithm is much faster with little impact on the performance.
  \item Several interesting phenomena are revealed. For example, compared with the state-of-the-art online algorithms, such as the greedy algorithm and MPC, the proposed algorithm can promote the energy sharing between charging stations and reduce the total cost by 8.12\%.
\end{enumerate}


\vspace{-0.cm}
\section{Mathematical Formulation}\label{sec:model}
\subsection{System Overview}
\begin{figure}[!htbp]
  \centering
  \vspace{-0.cm}  \includegraphics[width=0.4\textwidth]{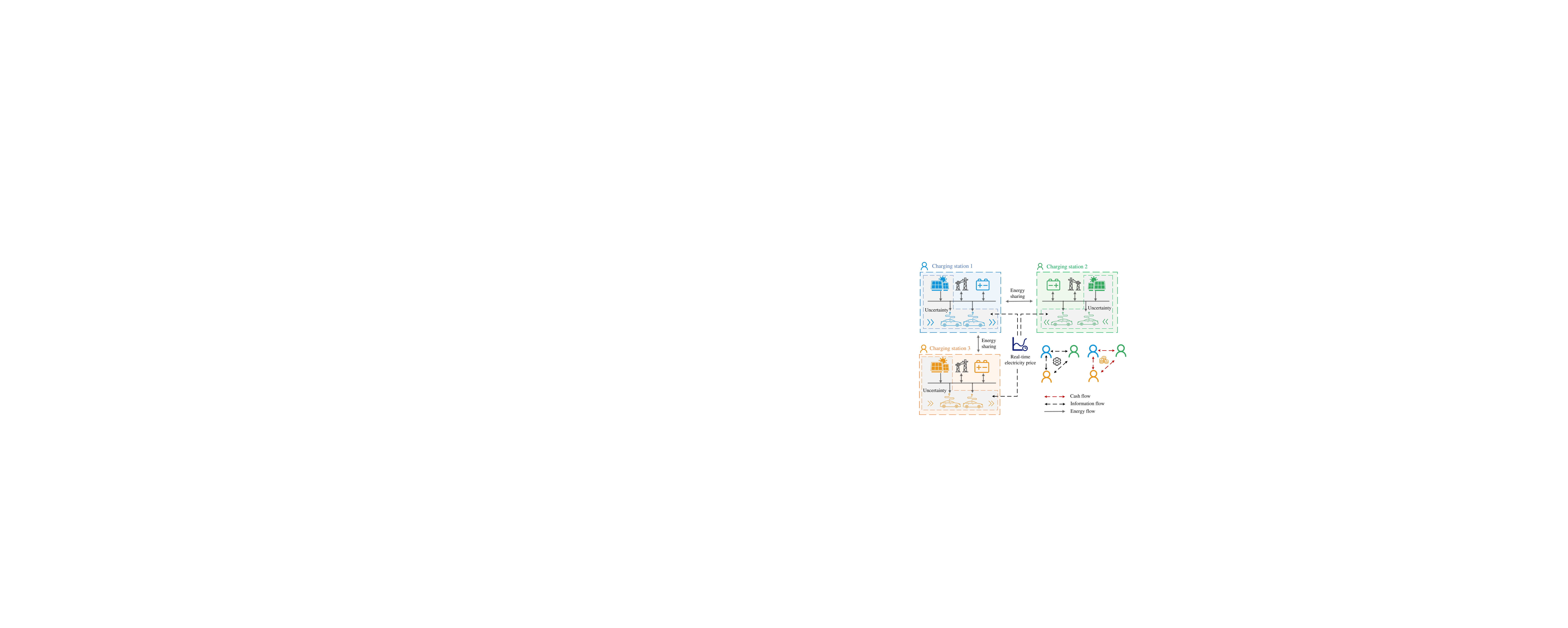}\\
  \vspace{-0.cm}
  \caption{Illustration of the energy flows between EV charging stations.}\label{fig:sysConf}
\end{figure}

We consider a microgrid/energy community with a group of EV charging stations working at discrete times $t\in\{1,...,T\}$, as shown in Fig.~\ref{fig:sysConf}.
Define $\mathcal{I}$ as the set of charging stations, each of which is indexed by $i\in\mathcal{I}$. For each charging station, its demand can be supplied by its own PV generation, battery storage, main grid, or the shared energy from neighboring charging stations.

\subsection{Modelling of Charging Station}
For each charging station, its modeling can be divided into two types:
1) Modeling of external energy flows, including energy trading with the main grid and energy sharing with neighboring stations;
2) Modeling of internal energy flows, including PV power generation, battery energy storage, and EV charging demand.
\subsubsection{Energy Trading with The Main Grid}
In addition to utilizing the energy from local PV, the charging station can purchase energy directly from the grid.
When local PV generation exceeds the charging demand, the charging station can sell the surplus energy to the grid.
We denote $p_{g,i}^b(t)$ and $p_{g,i}^s(t)$ as the purchasing and selling energy from/to the grid by charging station $i\in\mathcal{I}$ at time slot $t$.
Both $p_{g,i}^b(t)$ and $p_{g,i}^s(t)$ should satisfy physical limits:
\begin{align}
&  p_{g,i}^b(t) \ge 0,~ p_{g,i}^s(t) \ge 0, \label{equ:pgs}  
\end{align}

The energy trading cost of charging station $i\in\mathcal{I}$ at time slot $t$ incurred by trading with the main grid is:
\begin{equation}\label{equ:utilitycost}
C_{g,i}(t)=\left(p^b_{g,i}(t)c^b_{g}(t)-p^s_{g,i}(t)c^s_{g}(t)\right),
\end{equation}
where $c^b_{g}(t)$ and $c^s_{g}(t)$ are the unit electricity purchase and sale prices, respectively; and $c^{b,min}_{g}\leq c^b_{g}(t)\leq c^{b,max}_{g}$, $c^{s,min}_{g}\leq c^s_{g}(t)\leq c^{s,max}_{g}$.
In fact, electricity prices $c^b_{g}(t)$ and $c^s_{g}(t)$ are uncertain and only known at time slot $t$.
Further, they should meet the following requirement
\begin{equation}\label{inequ:price}
c^s_{g}(t)<c^b_{g}(t),
\end{equation}
which can ensure that charging stations won't arbitrage by buying from and selling back to the grid at the same interval.

\subsubsection{Energy Sharing with Other Charging Stations}
The renewable power supply and charging demand of adjacent charging stations may differ greatly due to the different charging patterns and facility capacities.
Therefore, charging station $i$ can share energy with its interconnected neighbor $j\in\mathcal{I}\backslash i$.
Through energy sharing, charging stations with surplus energy can transfer part of the energy to other charging stations that are lack of energy.
This can enhance the energy efficiency of the entire system and reduce the amount of energy bought/sold from/to the main grid.

Let $e_{i,j}(t)$ be the amount of energy that charging station~$i$ shares with charging station~$j$ at time slot $t$.
$e_{i,j}(t)>0$ means charging station~$i$ purchases energy from~$j$.
Conversely, if $e_{i,j}(t)<0$, it means that charging station~$i$ sells energy to~$j$.
The energy sharing between charging stations should meet the following coupling constraint:
\begin{equation}\label{equ:eij}
  e_{i,j}(t)+e_{j,i}(t)=0, \forall t\in\mathcal{T},\forall i\in\mathcal{I},\forall j\in\mathcal{I}\backslash i.
\end{equation}
In addition, $e_{i,j}(t)$ should be within the allowable range of minimum $e_{i,j}^{\min}$ and maximum $e_{i,j}^{\max}$:
\begin{equation}\label{equ:eijlbub}
  e_{i,j}^{\min}\leq e_{i,j}(t)\leq e_{i,j}^{\max}, \forall t\in\mathcal{T},\forall i\in\mathcal{I},\forall j\in\mathcal{I}\backslash i.
\end{equation}
To be fair, we have $e_{i,j}^{\min}=e_{j,i}^{\min}$ and $e_{i,j}^{\max}=e_{j,i}^{\max}$. The loss of energy during transmission is ignored as in \cite{wang}.
To encourage energy sharing among charging stations, we assume that the trading price $c_{sh}$ of energy sharing is lower than the electricity purchasing price from the grid ($c^b_{g}$) and greater than the selling price ($c^s_{g}$), i.e.,
\begin{equation}\label{inequ:sharePrice}
c^s_{g}<c_{sh}<c^b_{g}.
\end{equation}
Otherwise, charging stations have no incentive to participate in energy sharing.
The cost of station $i$ incurred by energy sharing at time slot $t$ is calculated by
\begin{equation}\label{equ:eshcost}
C_{sh,i}(t)=\sum\limits_{j\in\mathcal{I}\backslash i}e_{i,j}(t)c_{sh}(t).
\end{equation}

Budget balance is achieved for the overall charging stations system, i.e., $\sum_i C_{sh,i}(t)=0$, which is due to \eqref{equ:eij}.

\subsubsection{Local PV Power Generation}
Compared with other renewable energy sources such as wind power, PV generation is noiseless, and PV panels can be easily deployed on the rooftop of charging stations.
We let $p_{pv,i}(t)$ denote the generated PV energy of charging station $i\in\mathcal{I}$ in time slot $t$.
PV power generation is random and highly depends on the current solar irradiance.
We assume no prior knowledge of $p_{pv,i}$ nor its statistics.
In addition, unlike conventional generators, PV power generation has nearly zero cost.
\vspace{-0.cm}

\subsubsection{Battery Energy Storage}
Let $p^{d}_{b,i}(t)$ and $p^{c}_{b,i}(t)$ denote the battery storage discharging and charging energy of charging station $i\in\mathcal{I}$ at time slot~$t$.
$p^{d}_{b,i}(t)$ and $p^{c}_{b,i}(t)$ should meet the following physical constraints:
\begin{align}
  &0\leq p^{d}_{b,i}(t)\leq p_{b,i}^{d,\max},\label{equ:pbdch}\\
  &0\leq p^{c}_{b,i}(t)\leq p_{b,i}^{c,\max},\label{equ:pbcha}
\end{align}
where $p_{b,i}^{d,\max}$ and $p_{b,i}^{c,\max}$ are the maximum discharging and charging energy of battery, respectively.

The battery energy dynamics $E_{b,i}(t)$ can be represented by:
\begin{equation}\label{equ:Ebat}
E_{b,i}(t+1)=E_{b,i}(t)-\frac{1}{\eta_d}p^{d}_{b,i}(t)+\eta_{c}p^{c}_{b,i}(t),
\end{equation}
where $\eta_d$ and $\eta_c$ are the discharging and charging efficiency.
The battery energy should be always within its allowable range to avoid over-discharging or over-charging:
\begin{equation}\label{equ:EbatInequ}
E^{\min}_{b,i}\leq E_{b,i}(t)\leq E^{\max}_{b,i},
\end{equation}
where $E_{b,i}^{\min}$/$E_{b,i}^{\max}$ are the minimal/maximal level. 

We assume that the battery capacity is large enough to accommodate two consecutive charges or discharges, i.e.,

\noindent \textbf{A1}: $2\max\{p_{b,i}^{d,\max}/\eta_d,~ p_{b,i}^{c,\max}\eta_c\} \le E_{b,i}^{\max}-E_{b,i}^{\min},\forall i \in \mathcal{I}.$

Both frequent charging and discharging can cause battery degradation.
Here, the battery degradation cost is considered:
\begin{equation}\label{equ:batcost}
C_{b,i}( t )=c_{b,i}\left(p^d_{b,i}(t)+p^c_{b,i}(t)\right),
\end{equation}
where $c_{b,i}$ is the coefficient to measure the degradation caused by charging and discharging.

\subsubsection{EV Charging Demand}
EVs are flexible loads, whose charging demand can be partially shed in response to the power supply conditions.
For charging station $i\in\mathcal{I}$, its charging demand ${p_{d,i}}$ at time slot $t$ satisfies:
\begin{equation}\label{equ:pdBound}
p^{\min}_{d,i}\leq p_{d,i}(t)\leq p^{\max}_{d,i},
\end{equation}
where $p^{\min}_{d,i}$ is the minimum charging power at time slot $t$ required by EVs that cannot be shed, and $p^{\max}_{d,i}$ is the maximum charging energy preferred by EVs.

The charging demand shedding will cause discomfort for EVs, which is measured by
\begin{equation}\label{equ:costpd}
C_{l,i}(t)=\alpha_i(p^{\max}_{d,i}(t)-p_{d,i}(t))^2,
\end{equation}
where $\alpha_i$ is a positive coefficient to indicate the sensitivity of charging station $i$ to the unsatisfied charging demand $p^{\max}_{d,i}(t)-p_{d,i}(t)$.
In addition, to prevent excessive charging load shedding, an upper bound is imposed on the time-average charging load shedding rate:
\begin{equation}\label{equ:pdShedR}
r_{d,i}(t)=\frac{p^{\max}_{d,i}(t)-p_{d,i}(t)}{p^{\max}_{d,i}(t)-p^{\min}_{d,i}(t)},
\end{equation}
\begin{equation}\label{equ:pdAvg}
\lim\limits_{T\rightarrow\infty}\frac{1}{T}\sum\limits_{t=1}^{T}\mathbb{E}\left[r_{d,i}(t)\right]\leq \beta_{d,i},
\end{equation}
where $\beta_{d,i}$ is a positive constant, responsible for controlling the quality of charging services.

Each charging station $i\in\mathcal{I}$ should maintain energy balancing in each time slot $t$, which is described as
\begin{equation}\label{equ:pBala}
\begin{aligned}
p_{d,i}(t) = & p^b_{g,i}(t)- p^s_{g,i}(t) + e_{i,j}(t)\\
& + p_{pv,i}(t) + p^d_{b,i}(t) -p^c_{b,i}(t).
\end{aligned}
\end{equation}



\subsection{Offline Optimization of the Charging Stations System}\label{subsec:form}
Our goal is to minimize the long-term operational cost of the overall charging stations system, which can be achieved by co-optimizing the charging and discharging of battery storage, the energy exchange with the grid and adjacent stations, and the supply to charging demand over the scheduling time horizon.

Define the system state vector $\bf{s(t)}$ in time slot $t$ as a collection of PV generations, charging demands, and electricity prices, i.e.,
\begin{align}
\bf{s(t)} = \left(s_1(t),s_2(t),...,s_i(t)\right), \forall i\in\mathcal{I},\label{equ:state}\\
\mathbf{s_i(t)} = \left(p_{pv,i}(t),p^{\min}_{d,i}(t),p^{\max}_{d,i}(t), c_{g}(t)\right).
\end{align}
The system state is uncertain and we have no prior knowledge about their statistics.

Define the control decision variable vector $\mathbf{u(t)}$ as
\begin{align}\label{equ:ctrlU}
\bf{u(t)} =\left(u_1(t),u_2(t),...,u_i(t)\right), \forall i\in\mathcal{I},\\ \mathbf{u_i(t)}=\left(p^b_{g,i}(t),p^s_{g,i}(t), e_{ij}(t),p^d_{b,i}(t),p^c_{b}(t),p_{d,i}(t)\right).
\end{align}

The total operational cost at time $t$ includes the energy trading cost with grid, energy sharing cost with neighboring stations, battery degradation cost, and load shedding cost:
\begin{equation}\label{equ:cost}
C(t) = \sum\limits_{i\in\mathcal{I}}C_i(t)=\sum\limits_{i\in\mathcal{I}}\big(C_{g,i}(t)  + C_{sh,i}(t) + C_{b,i}(t)+ C_{l,i}(t)\big).
\end{equation}

Given the current system state $\mathbf{s(t)}$, the offline optimization of the charging stations system casts down to a stochastic optimization that minimizes the time-average total operational cost, as follows:
\begin{equation}\label{equ:P1}
\begin{aligned}
  \mathbf{P1}: \min \limits_{\mathbf{u(t)}} &\lim\limits_{T\rightarrow\infty}\frac{1}{T}\sum\limits_{t=1}^{T}\mathbb{E}\Big[C(t)\Big] \\
  \hbox{s.t.} & \: (\ref{equ:pgs}),(\ref{equ:eij}),(\ref{equ:eijlbub}),(\ref{equ:pbdch})-(\ref{equ:EbatInequ}),(\ref{equ:pdBound}),(\ref{equ:pdAvg}),(\ref{equ:pBala})
\end{aligned}
\end{equation}
where the expectation symbol $\mathbb{E}[\cdot]$ in the objective function is with regard to the random system state and corresponding control decisions.

The above offline optimization cannot be solved directly for two reasons:
1) It requires complete information of the future. However, in practice, the real-time electricity price, PV generation, and charging load are uncertain.
2) It runs in a centralized manner which requires private information of each charging station.

For the first issue, online algorithms are desired. The main difficulty that hinders online implementation is the time coupling constraints for battery energy dynamics (\ref{equ:Ebat}).
Common online algorithms, such as MPC and greedy algorithm, either require prior knowledge of uncertainties or simply ignore temporal coupling constraints.
To tackle the above challenges, this paper adopts the Lyapunov optimization framework.
Though Lyapunov optimization has been widely used for queueing problems, the time coupling constraints of battery storage makes it hard to apply in energy systems problems.
In Section~\ref{sec:form}, we will deal with this difficulty by virtual queues to reformulate the problem.
For the second issue, we develop an improved ADMM based distributed algorithm with \emph{iteration reduction} in Section~\ref{sec:distributed} to protect privacy and reduce computational burden. This enables fast execution of the algorithm and fits for real-time implementation.

\section{Online Algorithm}\label{sec:form}
In this section, we turn the offline problem \textbf{P1} into an online tractable form and prove that battery energy bound constraint (\ref{equ:EbatInequ}) is always satisfied though not explicitly considered in the online problem; the optimality gap with \textbf{P1} is also provided.
\subsection{Problem Modification}\label{subsec:pm}
Lyapunov optimization can solve a stochastic optimization problem with time-average constraints. Therefore, to solve $\mathbf{P1}$, 
the time-coupling constraint (\ref{equ:Ebat}) needs to be converted into a time-average one.
First, both sides of (\ref{equ:Ebat}) are summed over $t\in\{1,...,T\}$ and divide them by $T$ yields
\begin{equation}\label{equ:batEAvg1}
\frac{1}{T}\sum\limits_{t=1}^{T}\left[p_{b,i}(t)\right]= \frac{E_{b,i}(T+1)}{T}-\frac{E_{b,i}(1)}{T},
\end{equation}
\begin{equation}\label{equ:pb}
p_{b,i}(t)=-\frac{1}{\eta_d}p^{d}_{b,i}(t)+\eta_{c}p^{c}_{b,i}(t).
\end{equation}
We then take expectations on both sides of (\ref{equ:batEAvg1}) and take limits over $T$ to infinity yielding
\begin{equation}\small\label{equ:batEAvg2}
\lim\limits_{T\rightarrow\infty}\frac{1}{T}\sum\limits_{t=1}^{T}\mathbb{E}\left[p_{b,i}(t)\right]= \lim\limits_{T\rightarrow\infty}\mathbb{E}\left[\frac{E_{b,i}(T+1)}{T}\right]-\lim\limits_{T\rightarrow\infty}\mathbb{E}\left[\frac{E_{b,i}(1)}{T}\right],
\end{equation}
Due to the constraint (\ref{equ:EbatInequ}), both $E_{b,i}(T+1)$ and $E_{b,i}(1)$ are finite.
Thus, the right hand side of (\ref{equ:batEAvg2}) equal to zero, namely
\begin{equation}\label{equ:batEAvg3}
\lim\limits_{T\rightarrow\infty}\frac{1}{T}\sum\limits_{t=1}^{T}\mathbb{E}\left[p_{b,i}(t)\right]= 0.
\end{equation}
Now the problem $\mathbf{P1}$ turns to be
\begin{equation}\label{equ:P2}
\begin{aligned}
  \mathbf{P2}:\: \min \limits_{\mathbf{u(t)}} &\lim\limits_{T\rightarrow\infty}\frac{1}{T}\sum\limits_{t=1}^{T}\mathbb{E}\Big[C(t)\Big] \\
  \hbox{s.t.} & \: (\ref{equ:pgs}),(\ref{equ:eij}),(\ref{equ:eijlbub}),(\ref{equ:pbdch})-(\ref{equ:Ebat}),(\ref{equ:pdBound}),(\ref{equ:pdAvg}),(\ref{equ:pBala}),(\ref{equ:batEAvg3})
\end{aligned}
\end{equation}
In fact, constraint (\ref{equ:batEAvg3}) is a relaxed version of constraints (\ref{equ:Ebat})-(\ref{equ:EbatInequ}).
Thus, any feasible solution to $\mathbf{P1}$ is also a feasible solution to $\mathbf{P2}$, i.e., $\mathbf{P2}$ is less constrained than $\mathbf{P1}$.
The above relaxation step is critical which enables us to employ the Lyapunov optimization framework.

\subsection{Lyapunov Optimization Based Method}\label{subsec:lya}
Now the Lyapunov optimization framework can be used, the key ideas of which are briefly described below:
1) Construct virtual queues to transform the time-average constraints into queue stability constraints;
2) Define the Lyapunov function to obtain the Lyapunov drift and the drift-plus-penalty;
3) Minimize the upper bound of the drift-plus-penalty term.

\subsubsection{Construct Virtual Queues}
Lyapunov optimization aims to transform time-average constraints into queue stability constraints.
We define two virtual queues to deal with time-average constraints (\ref{equ:pdAvg}) and (\ref{equ:batEAvg3}) in problem $\mathbf{P2}$.

\emph{Battery Energy Queue:}
The virtual battery energy queue $B_i(t)$ is defined as follows:
\begin{equation}\label{equ:vB}
B_i(t)= E_i(t)-\theta_i(t),
\end{equation}
where $\theta_i(t)$ is a perturbation parameter designed to ensure the feasibility of constraint (\ref{equ:EbatInequ}), which will be explained later.
The dynamics of virtual battery energy queue is obtained as
\begin{equation}\label{equ:vBDyn}
B_i(t+1)= B_i(t)+p_{b,i}(t).
\end{equation}
Comparing (\ref{equ:vBDyn}) with (\ref{equ:Ebat}), it can be observed that $B_{i}(t)$ is actually a shifted version of $E_{b,i}(t)$ of battery storage.
But different from $E_{b,i}(t)$, the virtual energy queue $B_i(t)$ could be negative because of the perturbation parameter $\theta_i(t)$.
This shift can ensure that the constraint (\ref{equ:EbatInequ}) is met. 
In addition, due to (\ref{equ:batEAvg3}), it can be easily derived that the virtual battery energy queue $B_i(t)$ is also mean rate stable, i.e.,
\begin{equation}\label{equ:vBMrs}
\lim\limits_{t\rightarrow\infty}\frac{\mathbb{E}[B_i(t)]}{t}=0.
\end{equation}
It means that the queue does not grow faster than linearly with the time.

\emph{Charging Demand Shedding Queue:}
Similarly, we define a virtual queue $H_i(t)$ to deal with the time-average constraint in (\ref{equ:pdAvg}).
It evolves as follows:
\begin{equation}\label{equ:vDDyn}
H_{i}(t+1)= \max\{H_{i}(t)-\beta_{d,i},0\}+r_{d,i}(t).
\end{equation}
Let $H_{i}(0)=0$.
The virtual queue $H_{i}(t)$ accumulates the portion of unsatisfied charging load ratio at each time slot $t$.
Because of constraint (\ref{equ:pdAvg}), $H_{i}(t)$ is also mean rate stable.
\begin{equation}\label{equ:vDMrs}
\lim\limits_{t\rightarrow\infty}\frac{\mathbb{E}[H_{i}(t)]}{t}=0.
\end{equation}
Intuitively, to keep the virtual queue $H_{i}(t)$ stable, the arrival rate in time slot $t$, i.e., the shedding ratio $r_{d,i}(t)$, should not exceed the threshold $\beta_{d,i}$.

\subsubsection{Obtain Lyapunov Function and Drift-Plus-Penalty}
We define $\bf{\Theta(t)}=(\bf{B(t)}, \bf{H(t)})$ as the concatenated vector of virtual queues, where
\begin{align}
\mathbf{B(t)}=\left(B_1(t),...B_I(t)\right),\label{equ:vBvector}\\
\mathbf{H(t)}=\left(H_1(t),...H_I(t)\right).\label{equ:vD}
\end{align}
A Lyapunov function is then defined as follows
\begin{equation}\label{equ:LyaFun}
L(\mathbf{\Theta(t)})=\frac{1}{2}\sum\limits_{i\in{\mathcal{I}}}B_{i}(t)^2 + w\frac{1}{2}\sum\limits_{i\in{\mathcal{I}}}H_{i}(t)^2,
\end{equation}
where $w$ is the weight. $L(\mathbf{\Theta(t)})$ can be considered as a measure of the queue size.
A smaller value of $L(\mathbf{\Theta(t)})$ is preferred to push virtual queues $B_i(t)$ and $H_i(t)$ to be less congested.
Continually, the conditional one time slot Lyapunov drift is defined as follows:
\begin{equation}\label{equ:LyaDrift}
\Delta(\mathbf{\Theta(t)})=\mathbb{E}[L(\mathbf{\Theta(t+1)}) - L(\mathbf{\Theta(t)})|\mathbf{\Theta(t)}],
\end{equation}
The expectation is taken with respect to the random $\mathbf{\Theta(t)}$.

The Lyapunov drift is a measure of the expectation of the queue size growth given the current state $\mathbf{\Theta(t)}$.
Intuitively, by minimizing the Lyapunov drift, virtual queues are expected to be stabilized.
However, only minimizing the Lyapunov drift may lead to a high total operational cost.
Therefore, following the drift-plus-penalty method, we add the expected total operational cost (\ref{equ:cost}) in one time slot to (\ref{equ:LyaDrift}).
The drift-plus-penalty term is obtained as follows
\begin{equation}\label{equ:driftPlusP}
\Delta(\mathbf{\Theta(t)})+V\mathbb{E}[C(t)|\mathbf{\Theta(t)}],
\end{equation}
where $V$ is a weight parameter that controls the trade-off between virtual queues stability and operational cost minimization. We will show later in Proposition~\ref{prop-1} how to choose the value of this parameter.

\subsubsection{Minimizing the Upper Bound}
Problem \eqref{equ:driftPlusP} is still time-coupled due to the definition of $\Delta(\mathbf{\Theta(t)})$. To adapt to an online implementation,  instead of directly minimizing the drift-plus-penalty term, we minimize the upper bound to obtain the control decision. To be specific,
we first derive the one time slot Lyapunov drift expression as follows:

\begin{small}
\begin{align}
\small L(\mathbf{\Theta(t+1)}) -L(\mathbf{\Theta(t)}) =&\frac{1}{2}\sum\limits_{i\in{\mathcal{I}}}\Big\{\left[B_{i}(t+1)^2-B_{i}(t)^2\right]\nonumber \\
   &+ w\left[H_{i}(t+1)^2-H_{i}(t)^2\right]\Big\}.\label{equ:lyaD}
\end{align}
\end{small}
Based on the queue update equation (\ref{equ:vBDyn}) and (\ref{equ:vDDyn}), we have
\begin{align}
\frac{1}{2}&\left[B_i(t+1)^2-B_i(t)^2\right] \nonumber\\
&\leq B_{i}(t)p_{b,i}(t) + \frac{1}{2}\max\{(\eta_c p_{b,i}^{c,\max})^2,(\frac{1}{\eta_d}p_{b,i}^{d,\max})^2\},\label{equ:bUb}\\
\frac{1}{2}&\left[H_{i}(t+1)^2-H_{i}(t)^2\right] \nonumber\\
&\leq H_{i}(t)(r_{d,i}(t)-\beta_{d,i})
 + \frac{1}{2}(1+\beta_{d,i}^2).\label{equ:hUb}
\end{align}
Continually, we substitute inequalities (\ref{equ:bUb}) and (\ref{equ:hUb}) into drift-plus-penalty term yielding
\begin{small}
\begin{align}
&\Delta(\mathbf{\Theta(t)})+V\mathbb{E}[C(t)|\mathbf{\Theta(t)}] \leq A+\sum\limits_{i\in{\mathcal{I}}}B_{i}(t)\mathbb{E}\left[p_{b,i}(t)|\mathbf{\Theta(t)}\right]\nonumber\\
& + w\sum\limits_{i\in{\mathcal{I}}}H_{i}(t)\mathbb{E}\left[r_{d,i}(t)-\beta_{d,i}|\mathbf{\Theta(t)}\right]+V\mathbb{E}[C(t)|\mathbf{\Theta(t)}]\label{equ:dppInequ}
\end{align}
\end{small}
where $A=\frac{1}{2}\sum\limits_{i\in{\mathcal{I}}}\max\{(\eta_c p_{b,i}^{c,\max})^2,(\frac{1}{\eta_d}p_{b,i}^{d,\max})^2\}+\frac{1}{2}w\sum\limits_{i\in{\mathcal{I}}}(1+\beta_{d,i}^2)$ is a constant.

Recalling that the main principle of the Lyapunov optimization based method is to minimize the upper bound, we obtain the following real-time optimization problem
\begin{align}
  \mathbf{P3}:\: \min\: &\sum\limits_{i\in{\mathcal{I}}}B_{i}(t)p_{b,i}(t) +w\sum\limits_{i\in{\mathcal{I}}}H_{i}(t)r_{d,i}(t)+VC(t)\label{equ:P3}\\
  \hbox{s.t.} & \: (\ref{equ:pgs}),(\ref{equ:eij}),(\ref{equ:eijlbub}),(\ref{equ:pbdch})-(\ref{equ:pbcha}),(\ref{equ:pdBound}),(\ref{equ:pBala}),(\ref{equ:vBDyn}),(\ref{equ:vDDyn}).\nonumber
\end{align}
In each time slot $t$, given the current system state $\mathbf{s(t)}$ and virtual queue state $\mathbf{\Theta(t)}$, the proposed method determines the control decision $\mathbf{u(t)}$ by solving problem $\mathbf{P3}$. The term $wH_i(t)\beta_{d,i}$ is ignored since it is a constant in the problem. Hereafter, the original offline optimization problem $\bf{P1}$ has been decoupled into simple real-time (online) problems, which can be executed at each time slot without requiring a high-complex solver and a prior knowledge of uncertain states.

\subsection{Feasibility Guarantee and Performance Analysis}
Comparing constraints of $\mathbf{P1}$ with those of $\mathbf{P3}$, it can be observed that constraint (\ref{equ:EbatInequ}) is not considered in $\mathbf{P3}$. We may be concerned about whether the solution generated by $\mathbf{P3}$ is feasible for $\mathbf{P1}$.
In fact, by carefully designing the perturbation parameter $\theta_i(t)$, this bound constraint (\ref{equ:EbatInequ}) of battery energy state can be guaranteed, which is stated in the following proposition.

\begin{proposition} 
\label{prop-1}
When Assumption A1 holds, if we let
\begin{equation}\label{equ:theta}
\theta_i(t)=E^{\min}_{b,i}+\frac{1}{\eta_d}p^{d,\max}_{b,i}+\frac{V}{\eta_{c}}(c^{b,\max}_{g}+c_{b,i}), \forall t,
\end{equation}
where
\begin{align}
0\leq V \leq & V_{\max}\nonumber \\
               = & \min\limits_{i\in\mathcal{I}} \frac{E^{\max}_{b,i}-E^{\min}_{b,i}-\eta_{c}p^{c,\max}_{b,i}-\frac{1}{\eta_d}p^{d,\max}_{b,i}}
                      {\frac{1}{\eta_c}c^{b,\max}_g-\eta_{d}c^{b,\min}_{g}+c_{b,i}(\eta_d+\frac{1}{\eta_c})}\label{equ:Vmax},
\end{align}
the sequence of optimal solutions obtained by
the online problem \textbf{P3} satisfies constraint (\ref{equ:EbatInequ}).
\end{proposition}

The proof of Proposition \ref{prop-1} can be found in Appendix \ref{appendix-A}. Moreover, another important issue we care about is: what's the gap between the optimal solutions of online problem \textbf{P3} and offline problem \textbf{P1}? This is addressed below.
\begin{proposition}
\label{prop-2}
Denote the achieved long-term time-average cost objective value of $\mathbf{P1}$ and $\mathbf{P3}$ as $C^*$ and $\widehat{C}$, respectively. We have
\begin{equation}\label{equ:gap}
\widehat{C} - C^* \leq \frac{1}{V}A.
\end{equation}
where $A$ is a constant mentioned in (\ref{equ:dppInequ}).
\end{proposition}
The proof of Proposition \ref{prop-2} can be found in Appendix \ref{appendix-2}.
The optimality gap is affected by the control parameter $V$.
A bigger $V$ value can decrease the optimality gap but it increases the size of virtual queues.
In contrast, a smaller $V$ value makes queues more stable while leads to a larger optimality gap.

\section{Distributed Implementation}\label{sec:distributed}
Due to the need to protect privacy and reduce computational burden, we propose an ADMM based algorithm with iteration truncation to solve the energy sharing problem.

\subsection{Distributed Optimization Formulation}
First, auxiliary variable vectors $\boldsymbol{\varepsilon_{i}},\forall i \in \mathcal{I}$
are introduced, and each $\boldsymbol{\varepsilon_i}=\{\varepsilon_{i,j}(t),\forall t\in\mathcal{T},\forall j\in\mathcal{I}\backslash i\}$.
Then the constraint (\ref{equ:eij}) is replaced by
\begin{align}
  & \varepsilon_{i,j}(t)=e_{i,j}(t), \forall t\in\mathcal{T},\forall i\in\mathcal{I},\forall j\in\mathcal{I}\backslash i,\label{equ:eijnew}\\
  & \varepsilon_{i,j}(t)+\varepsilon_{j,i}(t)=0, \forall t\in\mathcal{T},\forall i\in\mathcal{I},\forall j\in\mathcal{I}\backslash i,\label{equ:eij+eji}
\end{align}

The corresponding augmented Lagrangian function of problem \textbf{P3} is written as follows:
\begin{align}\small\label{equ:lag}
  \mathcal{L}(\mathbf{u_i(t)}) =&\sum\limits_{i\in\mathcal{I}}\bigg[B_{i}(t)p_{b,i}(t)+wH_{i}(t)r_{d,i}(t)+VC_i(\mathbf{u_{i}(t)})\nonumber\\
  &+\sum\limits_{j\in\mathcal{I}\backslash i}\sum\limits_{t\in\mathcal{T}}\frac{\rho}{2}\bigg(e_{i,j}(t)-\varepsilon_{i,j}(t)+\frac{d_{i,j}}{\rho}\bigg)^2\bigg],
\end{align}
where $\rho$ is a positive penalty parameter of the augmented term, and $\boldsymbol{d}=\{d_{i,j}(t),\forall t\in\mathcal{T},\forall j\in\mathcal{I}\backslash i\}$ is the vector of associated dual variables of the equality constraint (\ref{equ:eijnew}).

Let $\mathcal{U}_i$ denote the constraints set of decision variable $\mathbf{u_i(t)}$ derived from (\ref{equ:P3}). After the above transformation, the problem can be solved in an iterative process.
The first step is each charging station $i$ updates its own variables $\mathbf{u_i(t)}$ by solving the local optimization problem:
\begin{align}\label{equ:step1}
\min\quad &B_{i}(t)p_{b,i}^{k+1}(t)+wH_{i}(t)r_{d,i}^{k+1}(t)+VC_i(\mathbf{u_{i}(t)}^{k+1})\nonumber\\
&+\sum\limits_{j\in\mathcal{I}\backslash i}\sum\limits_{t\in\mathcal{T}}\frac{\rho}{2}\bigg(e^{k+1}_{i,j}(t)-\varepsilon^{k}_{i,j}(t)+\frac{d^{k}_{i,j}(t)}{\rho}\bigg)^2,\\
\hbox{s.t.}\quad &\mathbf{u_{i}(t)}^{k+1}\in\mathcal{U}_i,\hbox{ given }\rho,\varepsilon_{i,j}^k(t),\nonumber
\end{align}
where $k$ is the iteration number.
After obtaining the $\boldsymbol{e_{i}}^{k+1}$ of each charging station, the second step is to update the auxiliary variables $\boldsymbol{\varepsilon_{i}}^{k+1},\forall i$ by solving the optimization problem:
\begin{align}\label{equ:step2}
\min&\quad \sum\limits_{i\in\mathcal{I}}\sum\limits_{j\in\mathcal{I}\backslash i}\sum\limits_{t\in\mathcal{T}}\frac{\rho}{2}\bigg(e^{k+1}_{i,j}(t)-\varepsilon_{i,j}(t)+\frac{d^{k}_{i,j}(t)}{\rho}\bigg)^2,\\
\hbox{s.t.}&\quad(\ref{equ:eij+eji}),\hbox{ given }\rho,e^{k+1}_{i,j}(t).\nonumber
\end{align}
By solving (\ref{equ:step2}), we can get the closed-form solution of $\boldsymbol{\varepsilon_{i}}^{k+1},\forall i \in \mathcal{I}$ as follows:
\begin{equation}\label{equ:step2-2}
\begin{aligned}
  \varepsilon^{k+1}_{i,j}(t)=&\frac{1}{2}\big((e^{k+1}_{i,j}(t)-e^{k+1}_{j,i}(t))
  +\frac{1}{\rho}(d^{k}_{i,j}(t)-d^{k}_{j,i}(t))\big),  \\
  \varepsilon^{k+1}_{j,i}(t)=&\frac{1}{2}\big((e^{k+1}_{j,i}(t)-e^{k+1}_{i,j}(t))
  +\frac{1}{\rho}(d^{k}_{j,i}(t)-d^{k}_{i,j}(t))\big).
\end{aligned}
\end{equation}
Based on  the updated $e^{k+1}_{i,j}(t)$ and $\varepsilon^{k+1}_{i,j}(t)$, the dual variable is updated in the third step by
\begin{align}
d^{k+1}_{i,j}(t)=d^{k}_{i,j}(t)+\rho\left(\varepsilon^{k+1}_{i,j}(t)-e^{k+1}_{i,j}(t)\right).\label{equ:step3}
\end{align}
The algorithm convergence criterion is set as:
\begin{equation}\label{equ:r}
  r=\|\boldsymbol{d}^{k+1}-\boldsymbol{d}^{k}\|\leq\delta.
\end{equation}
where $r$ is the error and $\delta$ is the accuracy tolerance.

\subsection{Efficiency Improvement for Real-time Implementation}\label{sec:efficient}
\begin{figure}[!htbp]
  \vspace{-0.cm}
  \centering
  \includegraphics[width=0.3\textwidth]{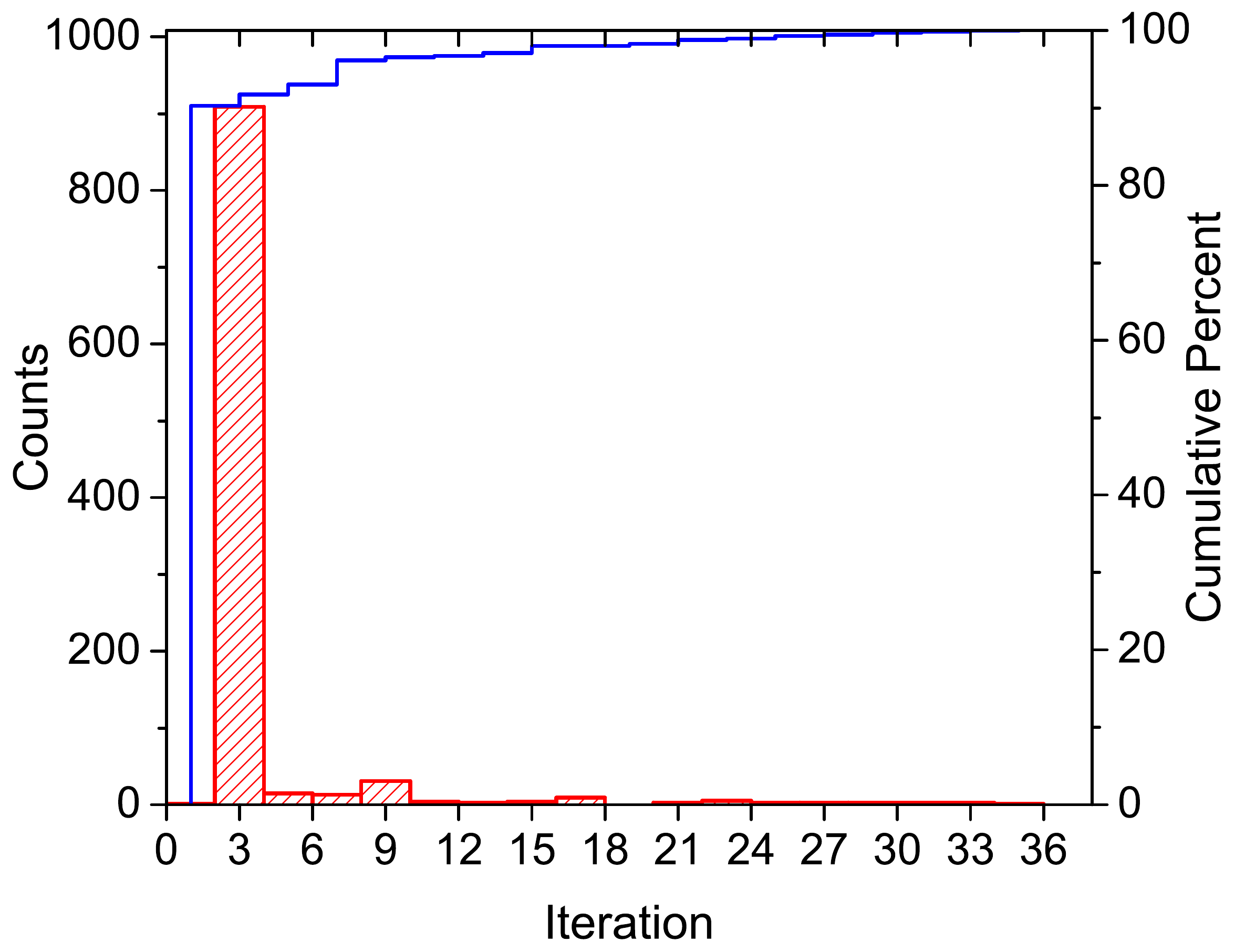}
  \caption{An example of the iteration numbers to reach convergence.}\label{fig:itrDstr}
\end{figure}

It may take a large number of iterations for the aforementioned algorithm to obtain the optimal strategy. However, in real-time application with a small time resolution, the algorithm may not have enough time to converge.
As shown in Fig.~\ref{fig:itrDstr}, the required iteration number are uncertain and depends on the situation in each time slot.
In the following, an iteration truncation method is developed to improve the efficiency of real-time energy sharing. It stops the algorithm when it exceeds a preset iteration threshold $k_{s}$.

Since the shared energy between a pair of charging stations $i$ and $j$ may not be balanced at the stop iteration $k_{s}$, here an energy balancing mechanism is proposed.
For the sake of secure operation, the shared energy quantities is determined by the minor one between $i$ and $j$.
For charging station $i$, if $|e^k_{i,j}|$ is greater than $|e^k_{j,i}|$, then charging station $i$ will decrease the shared energy $e^k_{i,j}$ according to the following update rule where $\hat{e}^k_{i,j}$ is the new shared energy after update,
\begin{equation}\label{equ:reBi}
  \hat{e}^k_{i,j} = - e^k_{j,i}, \forall i\in\mathcal{I},\forall j\in\mathcal{I}\backslash i,
\end{equation}
and charging station $j$ keeps no change,
\begin{equation}\label{equ:reBj}
  \hat{e}^k_{j,i} = e^k_{j,i}, \forall i\in\mathcal{I},\forall j\in\mathcal{I}\backslash i.
\end{equation}
Since $e_{i,j}^{\min}=e_{j,i}^{\min}$ and $e_{i,j}^{\max}=e_{j,i}^{\max}$, both $\hat e_{i,j}^k$ and $\hat e_{j,i}^k$ are feasible.
So far, the shared energy balance condition is guaranteed.
However, this shared energy update inevitably affect the system energy balance (\ref{equ:pBala}) obtained at iteration $k_{s}$ for charging station $i$.
The caused imbalance should be offset by adjusting dispatchable resources (grid or battery) of the charging station.
Here, we choose to adjust the purchasing/selling energy from/to the grid as follows
\begin{equation}\label{equ:reBpb}
  \left\{
  \begin{array}{ccc}
  p^b_{g,i} =p^b_{g,i} + (e^k_{i,j}+e^k_{j,i}) & \text{if }{ e^k_{i,j}\geq 0}\\ 
  p^s_{g,i} =p^s_{g,i}  - (e^k_{i,j}+e^k_{j,i}) & \text{if } {e^k_{i,j}<0}\\
  \end{array} \right.
\end{equation}

On the contrary, if $|e^k_{i,j}|$ is less than $|e^k_{j,i}|$, then charging station $i$ will maintain unchange, and charging station $j$ will update its shared energy to $\hat{e}^k_{j,i}$ according to the similar update rule (\ref{equ:reBi}) and adjust the grid power like (\ref{equ:reBpb}).

A complete description of the proposed distributed algorithm is shown in Algorithm \ref{algo:dc} that will be executed at each time slot.
The data exchanged among charging stations only include their shared energy while other private data (e.g. charging load, battery storage SOC) are well protected.

\begin{algorithm}[ht]
\caption{ADMM with Iteration Truncation}
\label{algo:dc}
\begin{algorithmic}[1]
\STATE Set iteration index $k=0$, convergence error tolerance $\delta>0$, penalty parameter $\rho>0$.
\STATE Initialize dual variables $\boldsymbol{d}^k=\boldsymbol{0}$, auxiliary variable $\{\boldsymbol{\varepsilon_{i}=0},i\in\mathcal{I}\}$
\REPEAT
\FOR{Each charging station $i\in\mathcal{I}$}
\STATE{Formulate the energy management subproblem for each charging station as (\ref{equ:step1}) and solves it with $\boldsymbol{\varepsilon}_{i}^k$ and $\boldsymbol{d}_i^k$}
\STATE Update $\boldsymbol{e_{i}}^{k+1}$
\ENDFOR
\STATE Update $\boldsymbol{\varepsilon}_i^{k+1}$ with $\boldsymbol{e}_{i}^{k+1}$ and $\boldsymbol{d}^k$ via (\ref{equ:step2-2})
\STATE Update dual variables $\boldsymbol{d}^{k+1}$ via (\ref{equ:step3}) with the new $\boldsymbol{\varepsilon}_i^{k+1}$ and $\boldsymbol{e_{i}}^{k+1}$.
\STATE Set $k=k+1$
\UNTIL convergence criterion (\ref{equ:r}) is satisfied or $k$ reaches the threshold $k_s$
\IF{$k_s$ is reached}
\STATE Update the shared energy according to (\ref{equ:reBi}) and (\ref{equ:reBj})
\STATE Adjust grid output based on (\ref{equ:reBpb}) to offset the imbalance caused by energy sharing update
\ENDIF
\STATE Update the battery energy state of each charging station
\STATE Update virtual energy queues (\ref{equ:vBDyn}) and load shedding queues (\ref{equ:vDDyn}) of each charging station
\end{algorithmic}
\end{algorithm}

\vspace{-0.cm}
\section{Simulation Results and Discussion}\label{sec:result}
In this section, we evaluate the effectiveness of the proposed method and compare it with the state-of-the-art approaches.
\vspace{-0.cm}
\subsection{Simulation Setup}
\begin{figure}[!htbp]
  \centering
  \includegraphics[width=0.35\textwidth]{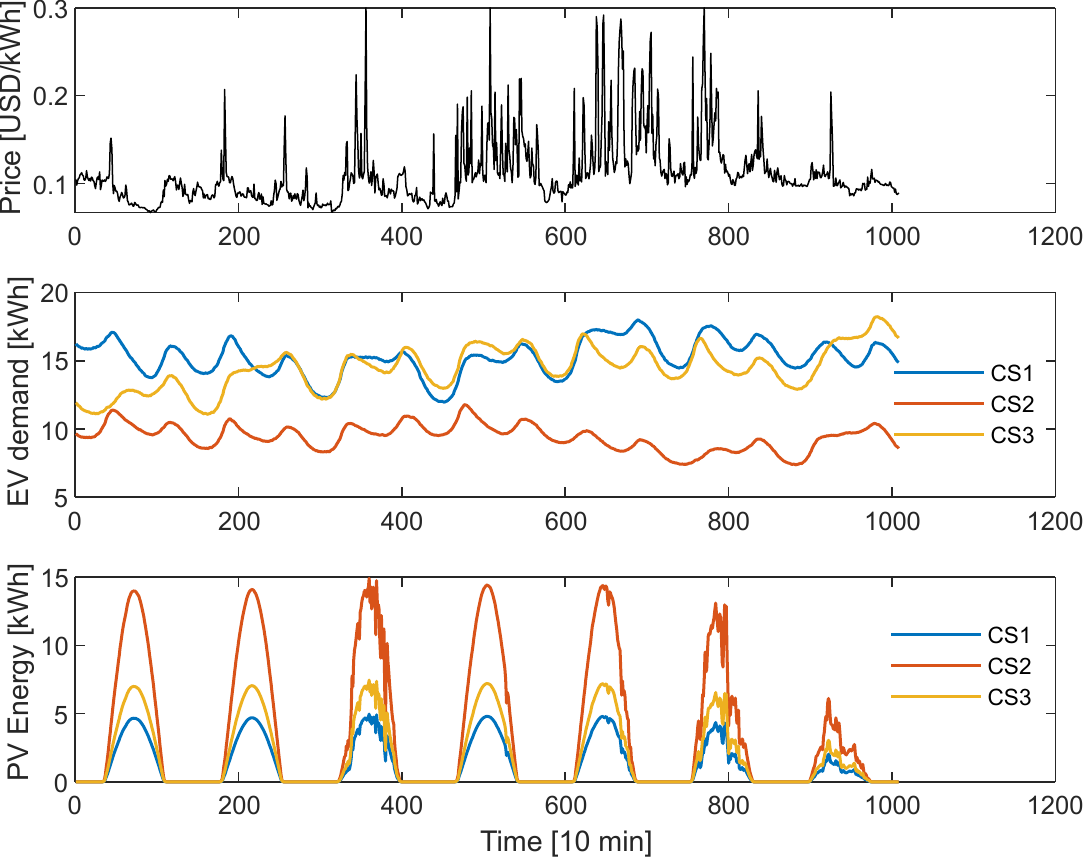}\\
  \vspace{-0.cm}
  \caption{Electricity price, EV charging demand, and PV energy output.}\label{fig:pvLoad}
\end{figure}

We first consider three charging stations (CS1, CS2, CS3) that interconnect with each other.
All of them are equipped with PV panels and battery storage, but their capacities are different.
The battery related parameters are $E_{b,1}=100\text{ kWh}, E_{b,2}=200\text{ kWh},E_{b,3}=200\text{ kWh},
p^{c,max}_{b,1}=p^{d,max}_{b,1}=10\text{ kW},
p^{c,max}_{b,2}=p^{d,max}_{b,2}=20\text{ kW},
p^{c,max}_{b,3}=p^{d,max}_{b,3}=20\text{ kW},
\eta_d=\eta_c=0.95, c_{b,1}=c_{b,2}=c_{b,3}=0.01.$
The used data of real-time electricity price~\cite{pjm} , PV power profile~\cite{nrel}, and EV charging load~\cite{pjm} is obtained from the real world to reflect the strong uncertainty.
They are shown in Fig.~\ref{fig:pvLoad}.
Specifically, we consider an entire period of 7 days (i.e., one week) with 10-minute per time slot, namely 1008 time slots in total.
The three charging stations have different charging demand patterns.
In particular, charging station 2 is with moderate charging demand while charging stations 1 and 3 are with heavy charging demand.
Charging station 2 has sufficient PV energy over the load while charging stations 1 and 3 need to import energy to meet their charging demand.
The feed-in tariff $c^s_{g}$ is set to be 0.01 \$/kWh at any time slot.
Referring to the commonly used mid-market rate mechanism, the internal energy sharing price $c_{sh}$ is chosen as an intermediate value between $c_g^s$ and $c_g^b$,
satisfying (\ref{inequ:sharePrice}).

\subsection{Simulation Results}
\subsubsection{Performance Comparison}
To show the advantage of the proposed distributed online algorithm, four widely used baselines in the literature are employed.
\begin{itemize}
  \item The first baseline (B1) individually operates each charging station while no energy sharing is allowed. Each station only minimizes the operational cost in the current time slot. Additionally, when the electricity price is below a threshold (0.1 USD), the battery storage will be charged at its maximum charging power.
  \item The second baseline (B2) uses the same greedy algorithm as B1, but energy sharing is incorporated. 
  \item The third baseline (B3) uses a MPC-based algorithm, which can look ahead and minimize the cost incurred in the prediction time window. Only the first step of the obtained control sequence will be applied. We assume B3 has accurate prediction over the future 6 time slots.
  \item The fourth baseline (B4) solves the offline optimization \textbf{P1}, assuming complete information of the future. Though not realistic, it provides a theoretical benchmark to verify the performance of other methods.
\end{itemize}

\begin{figure}[!htbp]
  \centering
  \includegraphics[width=0.35\textwidth]{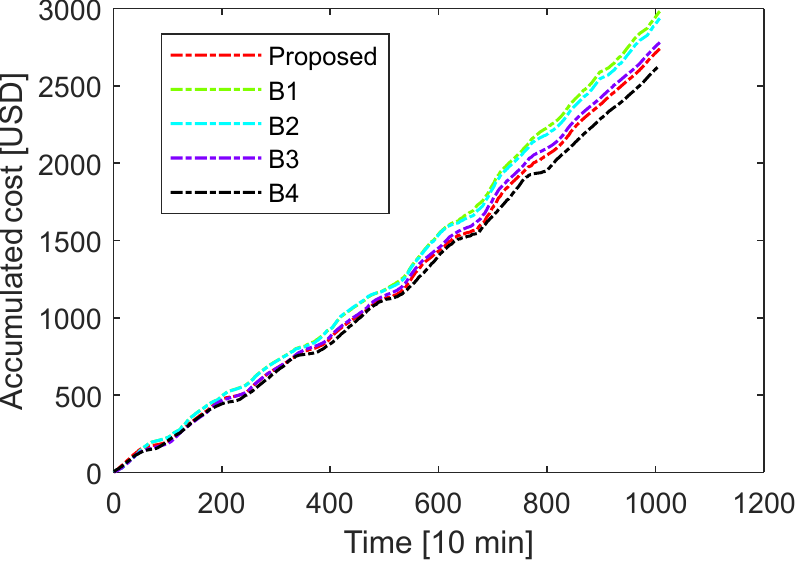}\\
  \vspace{-0.cm}
  \caption{Accumulated costs under five methods.}\label{fig:cmp}
\end{figure}

Fig.~\ref{fig:cmp} shows the accumulated operational costs over time under different methods.
The B1 algorithm has the worst performance and the highest cost.
The cost of B2 is lower than that of B1 thanks to energy sharing. With precise predictions for PV generation, charging load, and real-time energy pricing in the future 6 time slots, B3 outperforms B1 and B2.
However, the forecast and rolling procedure inevitably increases the computational complexity.
In contrast, our proposed algorithm performs better than B3 as time goes on and requires no future prediction, which is more practical.
In addition, the distance of the proposed algorithm from the greedy baseline B2 gradually widens over time.
This implies that the proposed algorithm has a long-term vision to achieve a better overall performance in the end. The offline optimization B4 has the best performance, which however is usually impossible in practice. 
Overall, the proposed algorithm can achieve a near offline optimal performance and is easy to implement.

TABLE~\ref{tab:costCmp} summarizes the costs under different methods.
The total cost of each charging station comprises the energy trading cost with the grid, battery cost, shedding cost and sharing cost.
Compared with B1, the proposed algorithm significantly reduces the total cost of three stations from 2981.75 USD to 2739.6 USD, with a larger drop of 8.12\%, while the reductions brought by B2 and B3 are merely 1.51\% and 6.74\%, respectively.
The result of the proposed algorithm is the closest to that of offline optimization B4.

\begin{table*}[!htbp]
  \centering
  \vspace{-0.cm}
    \caption{Cost comparison among B1, B2, B3, B4, and proposed algorithm (Unit: USD).}\label{tab:costCmp}
    \vspace{-0.cm}
   \begin{tabular}{@{}ccccccccc@{}}
     \hline\hline
                    & & Grid cost & Battery cost & Shedding cost & Sharing cost & Cost of each charging station & Total cost & Reduction\\
     \hline
     \multirow{3}{*}{B1}& CS1  & 1224.96 &14.96   & 94.50   & 0 & 1334.41 & \multirow{3}{*}{2981.75} & \multirow{3}{*}{-}\\
            &CS2  & 457.08   & 26.64  & 40.25   & 0 & 523.96   & \\
            &CS3  & 910.93   & 29.91  & 182.54 & 0 & 1123.38 &\\
     \hline
     \multirow{3}{*}{B2}&CS1 & 1010.08 & 14.96 & 94.50     & 193.84 & 1313.37 & \multirow{3}{*}{2936.59} &\multirow{3}{*}{1.51\%} \\
             &CS2 & 750.04 & 29.66   & 45.67   & -314.1  & 511.27  &  & \\
             &CS3 & 781.35 & 29.91   & 180.41 & 120.27 & 1111.95 &  & \\
    \hline
    \multirow{3}[0]{*}{B3} & CS1   & 976.30 & 6.47  & 94.50  & 206.27 & 1283.53 & \multirow{3}[0]{*}{2780.74} & \multirow{3}[0]{*}{6.74\%} \\
          & CS2   & 741.73 & 9.05  & 47.63 & -357.3 & 441.11 &       &  \\
          & CS3   & 708.86 & 10.72 & 185.49 & 151.03 & 1056.10 &       &  \\
    \hline
    \multirow{3}[0]{*}{B4} & CS1   & 553.3 & 13.63 & 94.5  & 568.73 & 1230.16 & \multirow{3}[0]{*}{2632.48} & \multirow{3}[0]{*}{11.71\%} \\
          & CS2   & 643.79 & 27.22 & 56.39 & -376.13 & 351.27 &       &  \\
          & CS3   & 1030.94 & 27.22 & 185.49 & -192.6 & 1051.04 &       &  \\
    \hline
     \multicolumn{1}{c}{\multirow{3}[0]{*}{Proposed}} & CS1   & 588.09 & 6.76  & 106.52 & 545.41 & 1246.78 & \multirow{3}[0]{*}{2739.60} & \multirow{3}[0]{*}{8.12\%} \\
          & CS2   & 954.40 & 14.11 & 62.83 & -611.02 & 420.32 &       &  \\
          & CS3   & 805.59 & 14.74 & 186.54 & 65.62 & 1072.50 &       &  \\
     \hline\hline
   \end{tabular}
\end{table*}

\subsubsection{Feasibility Analysis}

\begin{figure}[!htbp]
  \centering
  \includegraphics[width=0.35\textwidth]{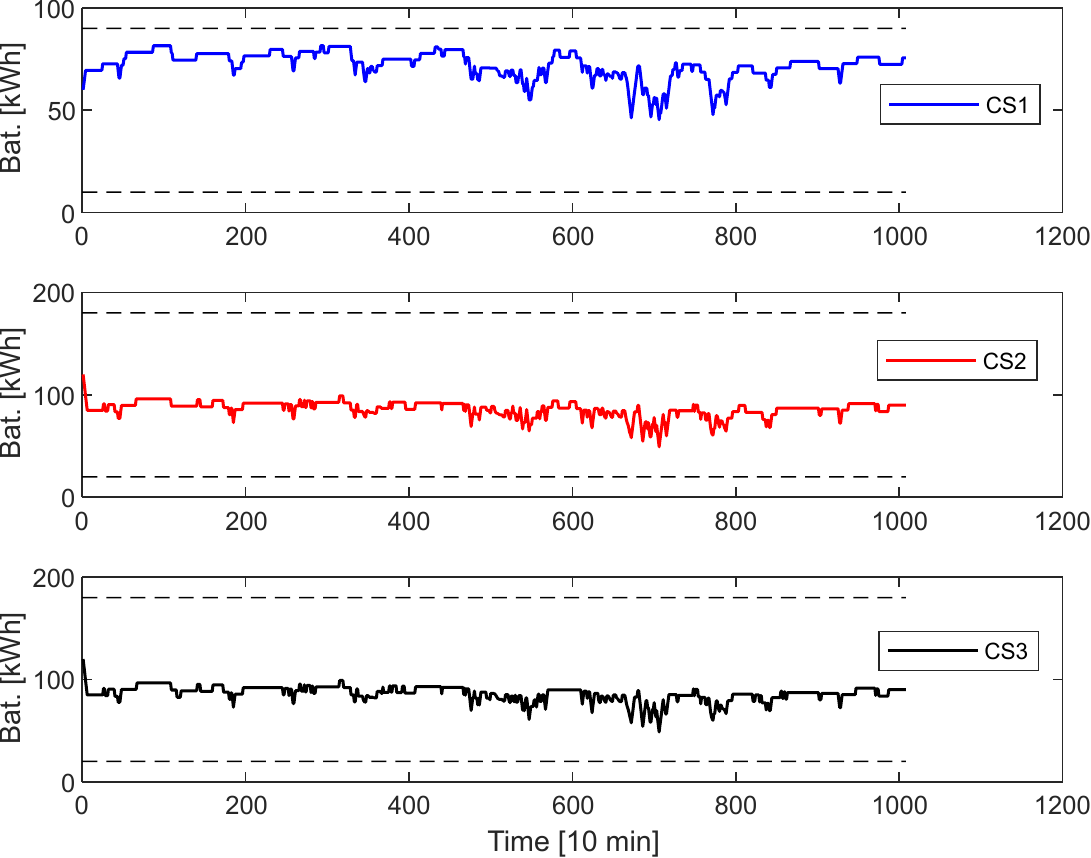}\\
  \vspace{-0.cm}
  \caption{Battery storage energy over time.}\label{fig:Ebat}
\end{figure}

We have proven that by appropriately designing the parameter $\theta_i(t)$, the battery energy bound constraint \eqref{equ:EbatInequ} can be met even though it is not explicitly considered in \textbf{P3}.
Here, simulation results are given to demonstrate the feasibility of the proposed algorithm.
Fig.~\ref{fig:Ebat} shows the energy evolution of battery storage over the time in each charging station.
As seen in the figure, all battery energy states are within the allowable range without exceeding the upper and lower bounds, which justifies Proposition \ref{prop-1}.
Meanwhile, the battery energy curve drops during the spikes of real-time electricity price. It means the battery discharges its stored energy when the price is very high and vice versa.
In addition, batteries remain at a medium energy level without large deviations, which enables it to quickly release or absorb energy when needed.
This is very helpful when working in future energy systems full of uncertainties.

\begin{figure}[!htbp]
  \centering
  \includegraphics[width=0.32\textwidth]{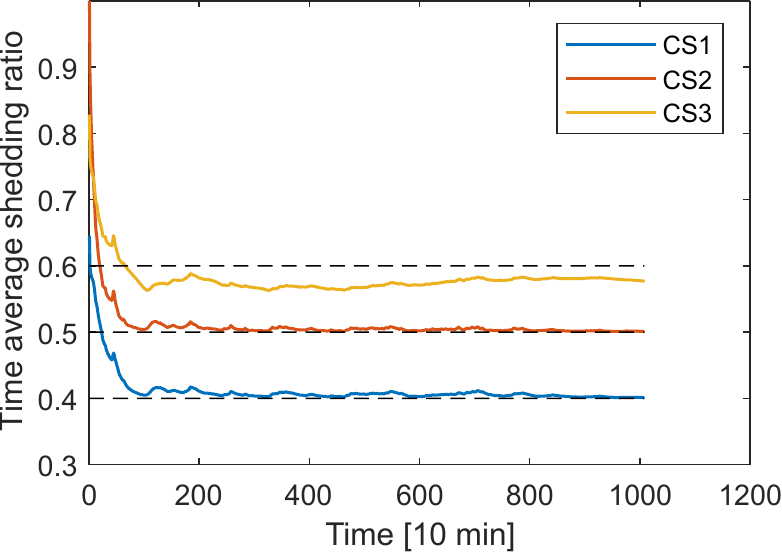}\\
  \vspace{-0.cm}
  \caption{Time-average charging load ratios.}\label{fig:loadRatio}
\end{figure}

Similarly, we replace the charging demand shedding ratio constraint \eqref{equ:pdAvg} by a virtual queue.
To examine its effectiveness, the time-average load shedding ratio over time is shown in Fig.~\ref{fig:loadRatio}.
At the beginning, the ratios are higher than the individual requirement of charging service quality.
Under the control of the proposed algorithm, they rapidly drop and meet the required time-average constraint as the time goes on.

\begin{figure}[!htbp]
  \vspace{-0.cm}
  \centering
    \subfigure[]{\label{fig:vload}\includegraphics[width=0.22\textwidth]{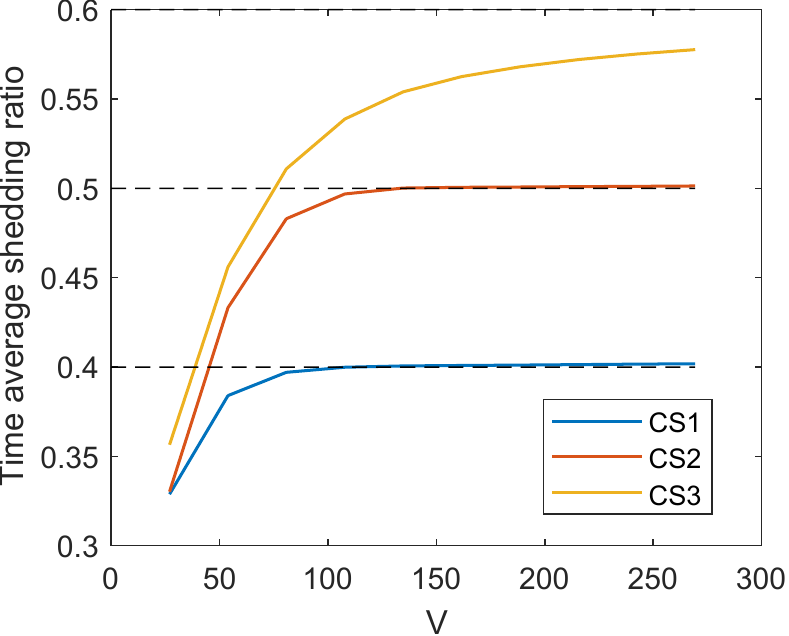}}
    \subfigure[]{\label{fig:vcost}\includegraphics[width=0.22\textwidth]{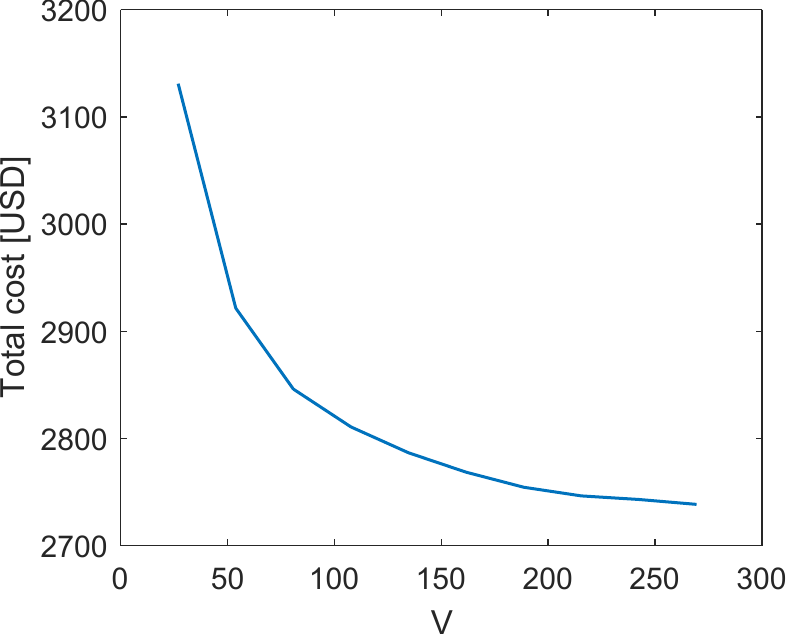}} \\
    \vspace{-0.cm}
  \caption{The impact of $V$ on the final time-average load shedding ratio and cost. (a) Final time-average load shedding ratio vs $V$. (b) Total cost vs $V$.}\label{fig:vimpact}
\end{figure}

Note that the above results are obtained under the weight parameter $V=V_{\max}$.
Recall that parameter $V$ controls the trade-off between stabilizing virtual queues and minimizing the total cost in the objective function~(\ref{equ:P3}).
Here we investigate the impact of $V$ by varying its value in the allowable range.
The results are shown in Fig.~\ref{fig:vimpact}.
The final time-average load shedding ratio and the total cost are both nonlinear function in $V$.
A larger $V$ would bring the time-average load shedding ratio closer to the requirement boundary, i.e., the queue turns to be more unstable, while the cost becomes lower.

\subsubsection{Energy Sharing Analysis}
\begin{figure}[!htbp]
  \vspace{-0.cm}
  \centering
  \includegraphics[width=0.4\textwidth]{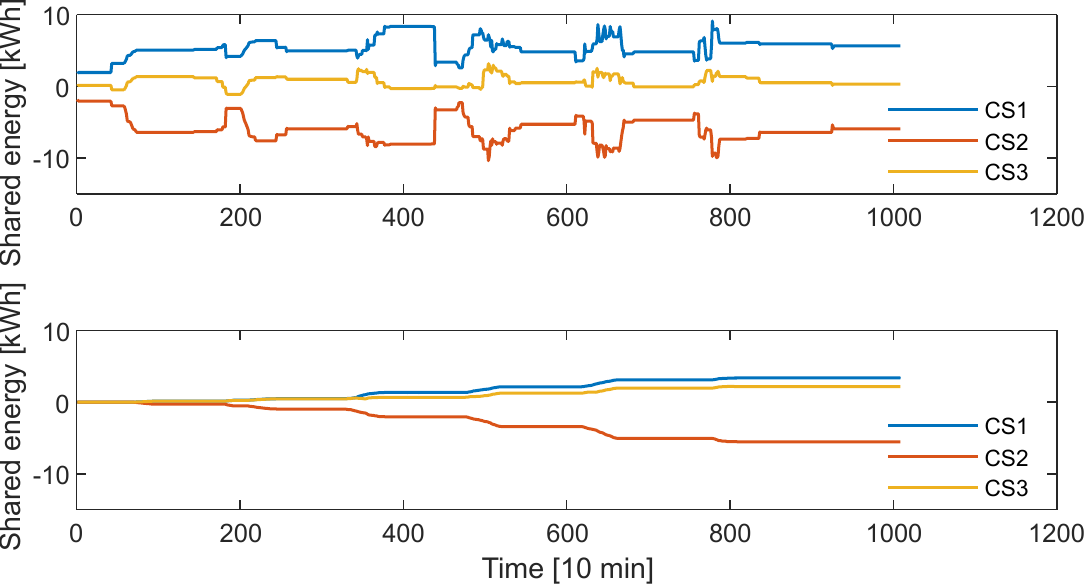}\\
  \vspace{-0.cm}
  \caption{Energy sharing result of the proposed algorithm (a) and B2 (b).}\label{fig:eShareCmp}
\end{figure}

The complementary nature of energy supply and demand between different charging stations makes energy sharing possible.
Fig.~\ref{fig:eShareCmp} (a) shows the energy sharing result of three charging stations in each time slot under the proposed algorithm.
Generally, CS1 and CS3 import energy from CS2, since CS2 has the most PV energy generation and the least charging load. 
Fig.~\ref{fig:eShareCmp} (b) shows the energy sharing solution under B2 algorithm.
Comparing the two figures, we can find that the proposed algorithm can greatly enhance the energy sharing between charging stations.

In addition, as seen in Table~\ref{tab:costCmp}, compared with B1 (no energy sharing case), B2 reduces the total cost by 1.51\% with the help of energy sharing.
However, the improvement in social welfare is marginal. The proposed algorithm further considers the long-term benefit through the Lyapunov optimization, achieving a significant cost reduction (8.12\%) without violating time coupling constraints.
Moreover, the proposed algorithm can reduce the cost of individual charging station compared with B1 and B2. Therefore, the charging stations have the incentive to participate in energy sharing.


\subsubsection{Convergence and Efficiency Analysis}
\begin{figure}[!htbp]
  \vspace{-0.cm}
  \centering
  \includegraphics[width=0.25\textwidth]{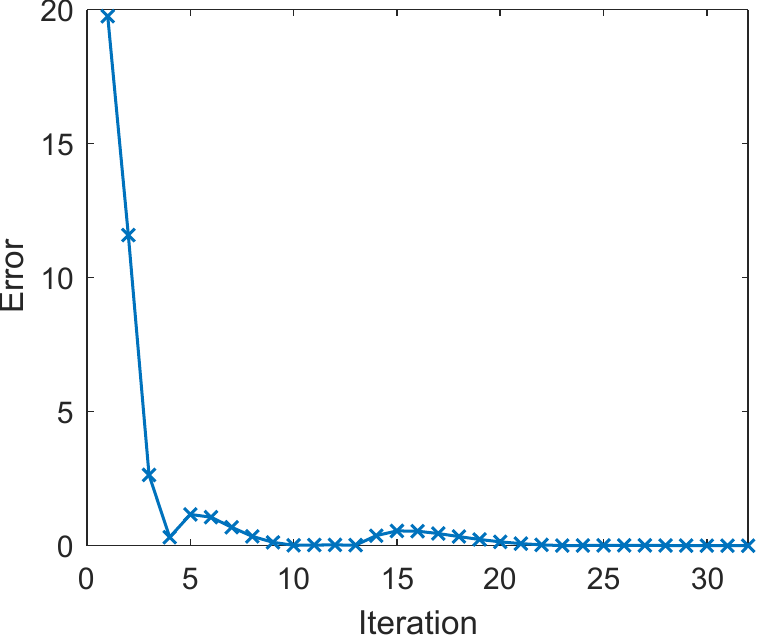}\\
  \vspace{-0.cm}
  \caption{An example of iterative convergence.}\label{fig:errCvg}
\end{figure}

Computational efficiency is a critical issue in real-time application.
Fig.~\ref{fig:itrDstr} has shown the iteration number to reach convergence without iteration truncation.
The corresponding error is given in Fig.~\ref{fig:errCvg}, which quickly drops at the beginning but then slowly converges. This allows us to stop the algorithm at a smaller number of iterations with little impact on the performance of the algorithm, which is the iteration truncation method in section \ref{sec:efficient}.
We set the threshold $k_s$ to 4, and compare the performance of the conventional ADMM algorithm with the proposed algorithm in terms of total cost, computing time and average iteration number.
As shown in Table~\ref{tab:effi}, the proposed algorithm achieves a higher computational efficiency (-26\%) and fewer iterations at slight cost increase (+0.12\%).

\begin{table}[!htbp]
  \centering
  \vspace{-0.cm}
    \caption{Performance comparison of A1 and A2 algorithm.}\label{tab:effi}
    \vspace{-0.cm}
   \begin{tabular}{@{}cccc@{}}
     \hline
           & Total cost (USD) & Time (s) & Average/Max. iteration \\
     \hline
     ADMM  & 2736.28                       & 115.40               & 3.01/34   \\
     Proposed  & 2739.60 (+0.12\%)     & 85.45  (-26\%)   & 2.20/4    \\
     \hline
   \end{tabular}
\end{table}

\subsubsection{Scalability}

\begin{figure}[!htbp]
  \vspace{-0.cm}
  \centering
    \subfigure[]{\label{fig:vsTime}\includegraphics[width=0.24\textwidth]{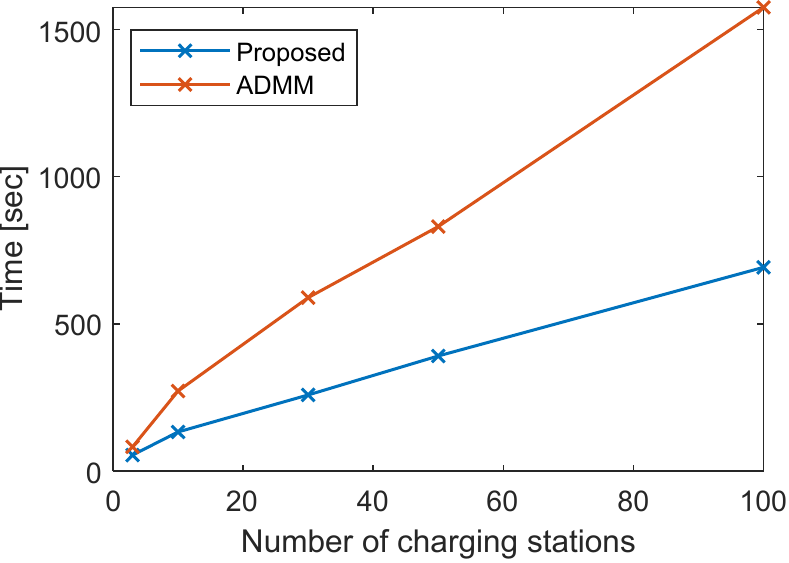}}
    \subfigure[]{\label{fig:vsCost}\includegraphics[width=0.23\textwidth]{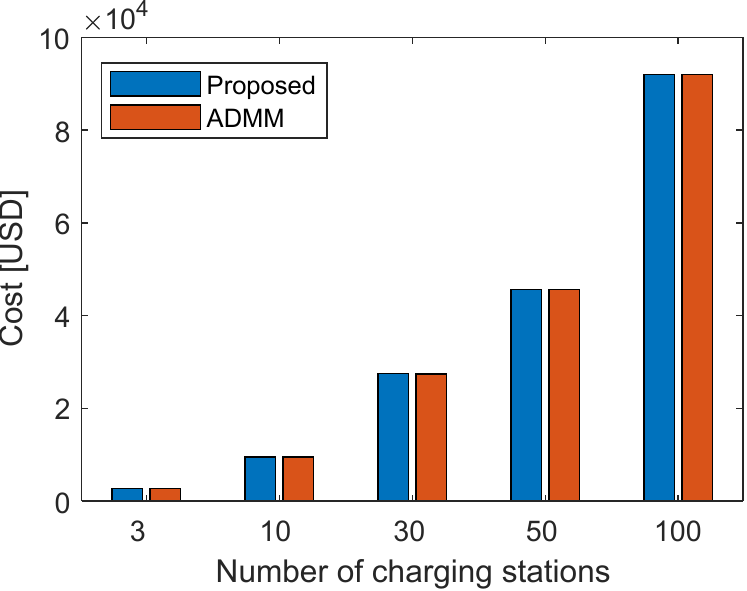}} \\
    \vspace{-0.cm}
  \caption{(a) Computational time. (b) Total cost.}\label{fig:scala}
\end{figure}

The above simulation results are obtained in the case of three charging stations.
In the following, we further investigate the scalability under massive charging stations. The computational time over the entire time slots (1008 time slots) and the total cost are used as criteria.
The results under different number of charging stations are shown in Fig.~\ref{fig:scala}.
As the number of charging stations increases, the proposed method takes much less time than the conventional ADMM with little impact on the total cost. 
For the extreme scenario with 100 charging stations, the proposed algorithm takes 692 s in total while the conventional ADMM needs 1574 s. The proposed algorithm is more scalable.


\section{Conclusion}\label{sec:conclu}
This paper proposes a distributed online algorithm to promote energy sharing among EV charging stations. Compared with the existing online algorithms, long-term benefits are considered through Lyapunov optimization technique where the time-coupling constraints are decoupled with the help of virtual queues. We provide guidance for selecting the parameters to ensure the satisfaction of battery related constraint. We also prove theoretically that the optimality gap between the proposed online algorithm and its offline counterpart is inversely proportional to the weight coefficient used in the drift-plus-penalty term. To protect privacy of individual charging stations, an improved ADMM algorithm with iteration truncation is proposed. Simulation results demonstrate the effectiveness and scalability of the proposed algorithm and have the following findings:

(1) The proposed online algorithm can achieve a nearly offline optimum, with a cost reduction of 6.7\% and 1.5\% compared to the greedy and MPC based algorithms, respectively

(2) Compared with the conventional ADMM algorithm, the computational time is reduced by 26\% with little sacrifice in total cost (+0.12\%).

(3) The participation of energy sharing can be enhanced by the proposed algorithm.
\vspace{-0.cm}

\appendix
\makeatletter
\@addtoreset{equation}{section}
\@addtoreset{theorem}{section}
\makeatother
\setcounter{equation}{0}  
\renewcommand{\theequation}{A.\arabic{equation}}
\renewcommand{\thetheorem}{A.\arabic{theorem}}

\subsection{Proof of Proposition \ref{prop-1}}
\label{appendix-A}
Suppose constraint (\ref{equ:EbatInequ}) holds at time slot $t$, we next prove that it also holds at time slot $t+1$ by induction.

\textbf{Case 1}: $E^{\min}_{b,i}\leq E_{b,i}(t) <\frac{1}{\eta_d}p^{d,\max}_{b,i}+E^{\min}_{b,i}$.
According to the energy balance constraint (\ref{equ:pBala}), we first derive the expression of $p^b_{g,i}$ and substitute it into the objective function of $\mathbf{P3}$ to eliminate the variable $p^b_{g,i}$.
Then, the partial derivative of objective function in $\mathbf{P3}$ with respect to $p^{c}_{b,i}(t)$ is
\begin{align*}
\frac{\partial P3(t)}{\partial p^{c}_{b,i}(t)}&=V\frac{\partial C(t)}{\partial p^{c}_{b,i}(t)}+B_i(t)\eta_c\\
                                                                     &\leq V(c^{b,max}_{g}+c_{b,i})+\left(E_{b,i}(t)-\theta_i(t)\right)\eta_c\\
                                                                     &=\left(E_{b,i}(t)-\frac{1}{\eta_d}p^{d,\max}_{b,i}-E^{\min}_{b,i}\right)\eta_c <0.
\end{align*}
Thus, the objective function is strictly decreasing with respect to $p^{c}_{b,i}(t)$.
Therefore, the optimal solution is $p^{c}_{b,i}(t)=p^{c,\max}_{b,i}$.
Since charging and discharging would not happen at the same time, we have $p^{d}_{b,i}(t)=0$.
Further, based on (\ref{equ:Ebat}), we have $E_{b,i}(t+1)= E_{b,i}(t)+\eta_c p^{c,\max}_{b,i}$ and hence
\begin{align*}
E^{\min}_{b,i}\leq E_{b,i}(t+1) \leq & E^{\min}_{b,i}+\frac{1}{\eta_d}p^{d,\max}_{b,i}+\eta_c p^{c,\max}_{b,i}\\
\leq & E^{\min}_{b,i}+E^{\max}_{b,i}-E^{\min}_{b,i}= E^{\max}_{b,i}.
\end{align*}
The third inequality is due to Assumption A1.

\textbf{Case 2}: $\frac{1}{\eta_d}p^{d,\max}_{b,i}+E^{\min}_{b,i}
\leq E_{b,i}(t) \leq V\left(\frac{1}{\eta_c}c^{b,\max}_g-\eta_{d}c^{b,\min}_{g}+c_{b,i}(\eta_d+\frac{1}{\eta_c})\right)+\frac{1}{\eta_d}p^{d,\max}_{b,i}+E^{\min}_{b,i}$.
Due to $V \leq V_{max} \leq \frac{E^{\max}_{b,i}-E^{\min}_{b,i}-\eta_{c}p^{c,\max}_{b,i}-\frac{1}{\eta_d}p^{d,\max}_{b,i}}
       {(1/\eta_c)c^{b,\max}_g-\eta_{d}c^{b,\min}_{g}+c_{b,i}(\eta_d+\frac{1}{\eta_c})}$,
we have
$E_{b,i}(t) \leq E^{\max}_{b,i}-\eta_{c}p^{c,\max}_{b,i}$.
Thus, based on the update (\ref{equ:Ebat}), we have
\begin{align*}
  E_{b,i}(t+1) & \leq E^{\max}_{b,i}-\eta_{c}p^{c,\max}_{b,i} -\frac{1}{\eta_d}p^{d}_{b,i}(t)+\eta_{c}p^{c}_{b,i}(t)\\
                       & \leq E^{\max}_{b,i}
\end{align*}
In addition, since $\frac{1}{\eta_d}p^{d,\max}_{b,i}+E^{\min}_{b,i}\leq E_{b,i}(t)$, we can obtain
\begin{align*}
  E_{b,i}(t+1) & \geq \frac{1}{\eta_d}p^{d,\max}_{b,i}+E^{\min}_{b,i}-\frac{1}{\eta_d}p^{d}_{b,i}(t)+\eta_{c}p^{c}_{b,i}(t)\\
                       & \geq E^{\min}_{b,i}
\end{align*}

\textbf{Case 3}: $V\left(\frac{1}{\eta_c}c^{b,\max}_g-\eta_{d}c^{b,\min}_{g}+c_{b,i}(\eta_d+\frac{1}{\eta_c})\right)+E^{\min}_{b,i}+\frac{1}{\eta_d}p^{d,\max}_{b,i} <E_{b,i}(t) \leq E^{\max}_{b,i}$.
Due to (\ref{equ:Vmax}), we have
$V\left(\frac{1}{\eta_c}c^{b,\max}_g-\eta_{d}c^{b,\min}_{g}+c_{b,i}(\eta_d+\frac{1}{\eta_c})\right)
+E^{\min}_{b,i}+\frac{1}{\eta_d}p^{d,\max}_{b,i}
\leq E^{\max}_{b,i}-\eta_{c}p^{c,\max}_{b,i}
< E^{\max}_{b,i}$.
Similar to \textbf{Case 1}, we then derive the partial derivative of the objective function of $\mathbf{P3}$ with respect to $p^{d}_{b,i}(t)$, i.e.,
\begin{align*}
&\frac{\partial P3(t)}{\partial p^{d}_{b,i}(t)}=    V\frac{\partial C(t)}{\partial p^{d}_{b,i}(t)}-B_i(t)\frac{1}{\eta_d}\\
&\leq V(c_{b,i}-c^{b,min}_{g})\\
&\quad-\left(E_{b,i}(t)-E^{\min}_{b,i}-\frac{1}{\eta_d}p^{d,\max}_{b,i}-\frac{V}{\eta_{c}}(c^{b,\max}_{g}+c_{b,i})\right)\frac{1}{\eta_d}\\
&=\Big[V\left(\frac{1}{\eta_c}c^{b,\max}_g-\eta_{d}c^{b,\min}_{g}+c_{b,i}(\eta_d+\frac{1}{\eta_c})\right)\\
&\quad+E^{\min}_{b,i}+\frac{1}{\eta_d}p^{d,\max}_{b,i} -E_{b,i}(t)\Big]\frac{1}{\eta_d}<0
\end{align*}
Thus, the objective function is strictly decreasing with respect to $p^{d}_{b,i}(t)$.
Therefore, the optimal solution is $p^{d}_{b,i}(t)=p^{d,\max}_{b,i}$.
Since charging and discharging cannot happen at the same time, we have $p^{c}_{b,i}(t)=0$.
According to (\ref{equ:Ebat}), we have $E_{b,i}(t+1)= E_{b,i}(t)-\frac{1}{\eta_d} p^{d,\max}_{b,i}$ and hence
\begin{align*}
E^{\min}_{b,i}
\leq E_{b,i}(t+1)
\leq & E^{\max}_{b,i}-\frac{1}{\eta_d} p^{d,\max}_{b,i}\leq E^{\max}_{b,i}
\end{align*}
Therefore, we have proved that the hard constraint (\ref{equ:EbatInequ}) still holds for all time slots. $\hfill\qedsymbol$

\setcounter{equation}{0}  
\renewcommand{\theequation}{B.\arabic{equation}}
\renewcommand{\thetheorem}{B.\arabic{theorem}}
\subsection{Proof of Proposition \ref{prop-2}}
\label{appendix-2}
Denote $\widehat{p_{b,i}}(t)$, $\widehat{r_{d,i}}(t)$ and $\widehat{C}(t)$ as the optimal results based on the optimal solution of $\mathbf{P3}$ in time slot $t$.
Denote $p^*_{b,i}(t)$, $r^*_{d,i}(t)$ and $C^*(t)$ as the optimal results of $\mathbf{P1}$ in time slot $t$.
According to (\ref{equ:dppInequ}), we have
\begin{align*}
&\Delta(\mathbf{\Theta}(t))+V\mathbb{E}[\widehat{C}(t)|\mathbf{\Theta}(t)] \leq A+\sum\limits_{i\in{\mathcal{I}}}B_{i}(t)\mathbb{E}\left[\widehat{p_{b,i}}(t)|\mathbf{\Theta}(t)\right]\nonumber\\
& + w\sum\limits_{i\in{\mathcal{I}}}H_{i}(t)\mathbb{E}\left[\widehat{r_{d,i}}(t)-\beta_i|\mathbf{\Theta}(t)\right]+V\mathbb{E}[\widehat{C}(t)|\mathbf{\Theta}(t)]\nonumber\\
&\leq
A+\sum\limits_{i\in{\mathcal{I}}}B_{i}(t)\mathbb{E}\left[p^*_{b,i}(t)|\mathbf{\Theta}(t)\right]\nonumber\\
& + w\sum\limits_{i\in{\mathcal{I}}}H_{i}(t)\mathbb{E}\left[r^*_{d,i}(t)-\beta_i|\mathbf{\Theta}(t)\right]+V\mathbb{E}[C^*(t)|\mathbf{\Theta}(t)]\\
&\leq
A+\sum\limits_{i\in{\mathcal{I}}}B_{i}(t)\mathbb{E}\left[p^*_{b,i}(t)\right]\nonumber\\
& + w\sum\limits_{i\in{\mathcal{I}}}H_{i}(t)\mathbb{E}\left[r^*_{d,i}(t)-\beta_i\right]+V\mathbb{E}[C^*(t)]
\end{align*}
Since the system states $\bf{s(t)}$ is i.i.d., $p^*_{b,i}(t)$ and $r^*_{d,i}(t)$ are also i.i.d. stochastic process.
Then, according to the strong law of large numbers, we obtain
\begin{align*}
&  \mathbb{E}[L(\mathbf{\Theta}(t+1)) - L(\mathbf{\Theta}(t))|\mathbf{\Theta}(t)]+V\mathbb{E}[\widehat{C}(t)|\mathbf{\Theta}(t)]\\
&  \leq A+\sum\limits_{i\in{\mathcal{I}}}B_{i}(t)\lim\limits_{T\rightarrow\infty}\frac{1}{T}\sum\limits_{t=1}^{T}\left\{p^*_{b,i}(t)\right\}\nonumber\\
& + w\sum\limits_{i\in{\mathcal{I}}}H_{i}(t)\lim\limits_{T\rightarrow\infty}\frac{1}{T}\sum\limits_{t=1}^{T}\left\{r^*_{d,i}(t)-\beta_i\right\}+V\mathbb{E}[C^*(t)]
\end{align*}
By taking expectation of the above inequality, we have
\begin{align*}
&  \mathbb{E}[L(\mathbf{\Theta}(t+1))] - \mathbb{E}[L(\mathbf{\Theta}(t))]+V\mathbb{E}(\widehat{C}(t))\\
&  \leq A+\sum\limits_{i\in{\mathcal{I}}}B_{i}(t)\lim\limits_{T\rightarrow\infty}\frac{1}{T}\sum\limits_{t=1}^{T}\mathbb{E}\left[p^*_{b,i}(t)\right]\nonumber\\
& + w\sum\limits_{i\in{\mathcal{I}}}H_{i}(t)\lim\limits_{T\rightarrow\infty}\frac{1}{T}\sum\limits_{t=1}^{T}\mathbb{E}\left[r^*_{d,i}(t)-\beta_i\right]+V\mathbb{E}[C^*(t)]\\
& \leq A+V\mathbb{E}[C^*(t)]
\end{align*}
By summing the above inequality over time slots $t\in\{1,2,\ldots,T\}$, we have
\begin{align*}
&  \sum\limits_{t=1}^{T}V\mathbb{E}[\widehat{C}(t)]\\
& \leq AT+V\sum\limits_{t=1}^{T}\mathbb{E}[C^*(t)]-\mathbb{E}[L(\mathbf{\Theta}(T+1))] + \mathbb{E}[L(\mathbf{\Theta}(1))].
\end{align*}
Since $L(\mathbf{\Theta}(T+1))$ and $L(\mathbf{\Theta}(1))$ are finite, we divide both sides of the above inequality by $VT$ and taking limits as $T\rightarrow\infty$ and obtain
\begin{align*}
\lim\limits_{T\rightarrow\infty}\frac{1}{T}\sum\limits_{t=1}^{T}\mathbb{E}(\widehat{C}(t))\leq \frac{A}{V}+\lim\limits_{T\rightarrow\infty}\frac{1}{T}\sum\limits_{t=1}^{T}\mathbb{E}(C^*(t))
\end{align*}
So far, We have finished the proof. $\hfill\qedsymbol$

\small{
\bibliographystyle{IEEEtran}
\bibliography{IEEEabrv,PaperRef}
}

\end{document}